\documentclass[12pt]{amsart}
\usepackage{amssymb}
\usepackage{amsthm}
\usepackage{amscd}

\usepackage{amsmath}
\usepackage{epic,eepic}
\usepackage{graphicx}
\DeclareMathAlphabet{\mathpzc}{OT1}{pzc}{m}{it}

\usepackage{mathrsfs}
\usepackage{latexsym}

\usepackage{pifont}

\theoremstyle{plain}
\newtheorem{lemma}{Lemma}[subsection]
\newtheorem{prop}[lemma]{Proposition}
\newtheorem{thm}[lemma]{Theorem}
\newtheorem{cor}[lemma]{Corollary}
\newtheorem{aplemma}{Lemma~A.\hspace{-1.5mm}}
\newtheorem{approp}{Proposition~A.\hspace{-1.5mm}}
\newtheorem{apthm}{Theorem~A.\hspace{-1.5mm}}
\newtheorem{apcor}{Corollary~A.\hspace{-1.5mm}}
\newtheorem{intthm}{Theorem}

\usepackage[all]{xy}
\usepackage {cancel}

\newcommand{\SSP}{\vspace{3mm}}
\newcommand{\LSP}{\vspace{5mm}}

\theoremstyle{definition}

\newtheorem{rema}[lemma]{Remark}

\newtheorem{remb}{Remark}

\newtheorem{defi}[lemma]{Definition}
\newtheorem{exa}[lemma]{Example}
\newtheorem{aprem}{Remark~A.\hspace{-1.5mm}}
\newtheorem{apdefi}{Definition~A.\hspace{-1.5mm}}
\newcommand{\bde}{\begin{defi}}
\newcommand{\ede}{\end{defi}\vspace{1mm}}
\newcommand{\ble}{\begin{lemma}}
\newcommand{\ele}{\end{lemma}}
\newcommand{\bpr}{\begin{prop}}
\newcommand{\epr}{\end{prop}}
\newcommand{\bt}{\begin{thm}}
\newcommand{\et}{\end{thm}}
\newcommand{\bco}{\begin{cor}}
\newcommand{\eco}{\end{cor}}
\newcommand{\bre}{\begin{rema}}
\newcommand{\ere}{\end{rema}}
\newcommand{\brea}{\begin{rema}}
\newcommand{\erea}{\end{rema}\vspace{1mm}}
\newcommand{\breb}{\begin{remb}}
\newcommand{\ereb}{\end{remb}\vspace{1mm}}
\newcommand{\bex}{\begin{exa}}
\newcommand{\eex}{\end{exa}}
\newcommand{\bpf}{\begin{proof}}
\newcommand{\epf}{\end{proof}\vspace{1mm}}

\newcommand{\bade}{\begin{apdefi}}
\newcommand{\eade}{\end{apdefi}}
\newcommand{\bale}{\begin{aplemma}}
\newcommand{\eale}{\end{aplemma}}
\newcommand{\bapr}{\begin{approp}}
\newcommand{\eapr}{\end{approp}}
\newcommand{\bat}{\begin{apthm}}
\newcommand{\eat}{\end{apthm}}
\newcommand{\baco}{\begin{apcor}}
\newcommand{\eaco}{\end{apcor}}
\newcommand{\bare}{\begin{aprem}}
\newcommand{\eare}{\end{aprem}}

\newcommand{\be}{\begin{enumerate}}
\newcommand{\ee}{\end{enumerate}}
\newcommand{\bcd}{\[\begin{CD}}
\newcommand{\ecd}{\end{CD}\]}
\newcommand{\bit}{\begin{itemize}}
\newcommand{\eit}{\end{itemize}}
\newcommand{\bq}{\begin{quote}}
\newcommand{\eq}{\end{quote}}
\newcommand{\ba}{\begin{array}}
\newcommand{\ea}{\end{array}}

\newcommand{\mcC}{\mathcal{C}}
\newcommand{\mcD}{\mathcal{D}}

\newcommand{\mcO}{\mathcal{O}}
\newcommand{\mcP}{\mathcal{P}}

\newcommand{\mcS}{\mathcal{S}}
\newcommand{\mcT}{\mathcal{T}}

\newcommand{\msU}{\mathscr{U}}

\newcommand{\migi}{\rightarrow}

\newcommand{\isom}{\stackrel{\sim}{\migi}}

\newcommand{\migiincl}{\hookrightarrow}

\newcommand{\mr}{\mathrm}
\newcommand{\hidden}[1]{\,}

\begin{document}

\title[On 
 the category of structure species]{
On the category of structure species}
\author{Yasuhiro Wakabayashi}
\date{}
\markboth{Yasuhiro Wakabayashi}{}
\maketitle
\footnotetext{Y. Wakabayashi: Department of Mathematics, Tokyo Institute of Technology, 2-12-1 Ookayama, Meguro-ku, Tokyo 152-8551, JAPAN;}
\footnotetext{e-mail: {\tt wkbysh@math.titech.ac.jp};}
\begin{abstract}

The purpose of the present paper is to make a mathematical study of the differences and relations among possible structures inherent in an object, as well as of the whole structure constituted by them (i.e., the structure of structures), against the background of the structuralism by Claude  L\'{e}vi-Strauss and others. Our discussion focuses on Blanchard's categorical reformulation of the notion of structure species introduced originally by Bourbaki. The main result  of the present  paper asserts that a category can be reconstructed, up to a certain slight indeterminacy, from  the category of structure species on it. This result is partially motivated by various reconstruction theorems that have been shown in the context of anabelian geometry.

\end{abstract}
\tableofcontents 

\section{Introduction: L\'{e}vi-Strauss' structuralism}

\LSP
\subsection{} \label{S01}

{\it Structuralism} is a general theory of culture and methodology
 that
focuses  on relationships rather than individual objects,  or alternatively, where objects are defined by the set of relationships of which they are part and not by the qualities possessed by them taken in isolation.
According to this mode of knowledge, 
 phenomena of human life are not intelligible except through their interrelations.
These relations constitute a {\it structure}, and behind local variations in the surface phenomena there are constant laws of abstract structure (cf.  ~\cite{Cal}, ~\cite{Rot}, ~\cite{SBla}).

Structuralism in Europe developed in the early 20th century
 from insights in  the field of linguistics of, mainly, 
  Ferdinand de Saussure and the subsequent Prague, Moscow,
and Copenhagen schools of linguistics. 
 After World War II, an array of scholars in the humanities borrowed Saussure's concepts for use in their respective fields. 
 By the early 1960s, structuralism as a movement was coming into its own and some believed that it offered a single unified approach to human life that would embrace all disciplines. 
 Claude L\'{e}vi-Strauss, a French anthropologist, was the first such scholar, sparking a widespread interest in structuralism.

L\'{e}vi-Strauss defined ``structure" as
a whole consisting of elements and relations between elements, which retain their invariant properties through a series of transformations.
To be precise, a model  with structural value satisfies several requirements described as follows (cf. ~\cite{Lev2}, Part 5,  Chap.\,XV):
 \begin{itemize}
 \item[(i)]
The structure exhibits the characteristics of a system. It is made up of several elements, none of which can undergo a change without effecting changes in all the other elements.
 \item[(ii)]
 For any given model there should be a possibility of ordering a series of transformations resulting in a group of models of the same type.
 \item[(iii)]
 The above properties make it possible to predict how the model will react if one or more of its elements are submitted to certain modifications.
 \item[(iv)]
 The model should be constituted so as to make immediately intelligible all the observed facts.
 \end{itemize}
 The structures listed, at least, as illustrations are kinships, political ideologies, mythologies, ritual, art, code of etiquette, and even cooking (cf. ~\cite{Rot}). The discernment of these structures and their comparative analysis, which takes into account their distribution  both historically and geographically, is indeed the subject matter of structural anthropology.
In ``{\it The Elementary Structures of Kinship}" (cf. ~\cite{Lev3}), Claude L\'{e}vi-Strauss examined kinship systems from a structural point of view and demonstrated how apparently different social organizations were different permutations of a few basic kinship structures.
The kinship system provides a way to order individuals according to certain rules; social organization is another way of ordering individuals and groups; social stratifications, whether economic or political, provide us with a third type; and all these orders can themselves be ordered by showing the kinds of relationship which exist among them, i.e.,  how they interact with one another on both the synchronic and the diachronic levels. 
One may construct models valid not only for one type of order (kinship, social organization, economic relations, etc.) but where numerous models for all types of order are themselves ordered inside a total model.
 In ~\cite{Rot}, Chap.\,II, N. Rotenstreich explained  that 
\begin{quotation}
{\it We do not start out with scattered concepts: we start with structures and move to further structures. We start with order
and move from another order or to an order of orders.
There is a structure to the relation between these orders or structures. This ``{\bf structure of structures}"
\footnote{L\'{e}vi-Strauss  developed  a similar idea with the term ``order of orders" (cf. ~\cite{Lev2}).}
 is not just a static relation of coexistence, i.e. language beside kinship, etc., or even
not one of subordination whereby a narrow structure such as rites is comprised in or is secondary to a wider structure such as society.}
\end{quotation}

\LSP
\subsection{} \label{S033}

In the present paper, we make a mathematical study of the differences and relations among possible structures inherent in an object, as well as  the whole structure constituted by them (i.e., the structure of structures), against the background of the structuralism by   L\'{e}vi-Strauss and others. 
Our discussion focuses on
 some objects treated in structural theory, starting with  Bourbaki.
As is well known,   Bourbaki, 
 the collective pseudonym of a group of predominantly French mathematicians,   undertook, in the mid-20th century, 
the task of making a unified development of central parts of
 modern mathematics in largely formalized language.\footnote{In 1943, Andr\'{e} Weil, one of the collaborators  in Bourbaki, met Claude L\'{e}vi-Strauss in New York, 
 which led to a small collaboration.
  By using a mathematical model based on group theory, Weil 
  described marriage rules for four classes of people within Australian aboriginal society.
  This contribution appeared in an appendix of L\'{e}vi-Strauss's book (cf. ~\cite{Lev3}).}
 This resulted in a long series of books  ``{\it \'{E}l\'{e}ments de Math\'{e}matique}" that became very influential.
In a manifesto written by  Bourbaki  in 1950, some main principles of their structuralist view of mathematics were presented.

In one of these books (cf. ~\cite{Bou}),
Bourbaki developed  their theory of
structures in a set-theoretic manner; the building blocks are called {\it structure species}.
Roughly speaking, a  structure species is a set, or a collection of sets, endowed with relations and operations not only among their members
but also among collections of elements of these sets, relations among them, etc. 
A basic example is the structure species of ordered sets,  where, from a set $S$, we obtain (by a suitable echelon construction scheme) the power set $P (S \times S)$ of the product   $S \times S$ and  a binary relation $s \in P (S \times S)$ (called the typical characterization) equipped with a relation 
 describing the axiom of an ordered set.

In ~\cite{Bla}, Blanchard introduced the concept of structure species on a category very close to the concept of structure species in the sense of Bourbaki and proved that 
it  is equivalent to the concept of structure species in the sense of Sonner (cf. ~\cite{Son}).
The purpose of the present paper is  to consider   ``structure of structures" formulated  in terms of Blanchard's structure species and investigate how this concept can capture the essence of things.
 That is, we discuss the issue of how much information concerning a given category is contained in the knowledge about the {\it structure of structure species} on that category.

The following is the main result of the present paper, which 
asserts that a category can be reconstructed, up to a certain slight indeterminacy, 
 from the category of structure species on it (the proof will be given in  \S\,\ref{S44}):

\SSP
\begin{intthm} \label{E223}
Let us fix a universe  $\msU$.
Let $X$ and $X'$ be connected $\msU$-small categories.
Denote by $\mcS p_X$ and $\mcS p_{X'}$ the categories of structure species on $X$ and $X'$ respectively.
Also, denote by $X^{\mr{op}}$ the opposite category of $X$.
Then, the following  two conditions are equivalent to each other:
\begin{itemize}
\item[(a)]
$X \cong X'$ or $X^\mr{op} \cong X'$.
\item[(b)]
$\mcS p_{X} \cong \mcS p_{X'}$.
\end{itemize}
Here, given two categories $\mcC$, $\mcD$, we write $\mcC \cong \mcD$ if $\mcC$ is equivalent to $\mcD$.
 \end{intthm}

\LSP
\subsection{} \label{S037}

Finally, we remark  that
Theorem  \ref{E223} can be regarded as a variant of various reconstruction theorems that have been shown in the context of anabelian geometry  (cf. Remark \ref{R450}).
In 1983, A. Grothendieck wrote a letter to G. Faltings (cf. ~\cite{Gro}), outlining what is today known as the anabelian conjectures. 
(Many of the claims based on these conjectures have now been proved by mathematicians.)
These conjectures concern the possibility of reconstructing  certain arithmetic varieties (e.g., hyperbolic curve over a number or a $p$-adic local field) from their \'{e}tale fundamental groups.
Our study 
is partially motivated by the anabelian philosophy of Grothendieck;
this is because a structure species, or equivalently a constructive functor (cf. Definition \ref{D02}, (i)), on a category $X$  may be regarded, in some sense, as categorical realization of coverings  over $X$ (cf. ~\cite{Son}, \S\,3).
Moreover,  if $X$ is a groupoid, then one can interpret Theorem \ref{E223}  as a reconstruction assertion for  $X$ by means of  the fundamental group associated with the Galois category of constructive functors over  $X$ (cf. Remark \ref{Err0}).
Similar category-theoretic reconstructions can be found in ~\cite{Moc1}, ~\cite{Moc2}, and ~\cite{Wak}.

\LSP
\subsection*{Acknowledgements} \label{S0341}
The author was partially supported by Grant-in-Aid for Scientific Research (KAKENHI No. 18K13385, 21K13770).

\vspace{10mm}
\section{Structures species and constructive functors} 
\SSP

In this section, we recall structure species on a category defined by Blanchard and the  equivalence  of structure species and constructive functors.
After that, we examine constructive functors on a groupoid.

\LSP
\subsection{Preliminaries on categories} \label{SS01}
Let $\mcC$ be a category.
We denote by $\mr{Ob}(\mcC)$ (resp., $\overline{\mr{Ob}}(\mcC)$; resp., $\mr{Mor}(\mcC)$) the set of objects (resp., the set of isomorphism classes of objects; resp., the set of morphisms) in $\mcC$.
For two objects $a, b \in \mr{Ob}(\mcC)$, we denote by $\mr{Mor}_\mcC (a, b)$, or $\mr{Mor} (a, b)$,  the set of morphisms $a \migi b$ in $\mcC$.
Also, we  write
\begin{align}
\mr{Mor}\{ a, b\} := \mr{Mor} (a, b) \cup \mr{Mor} (b, a)
\end{align}
(hence  $\mr{Mor}\{a, b\} = \mr{Mor}\{b, a\}$ and   $\mr{Mor}\{a, a\} = \mr{Mor}(a, a)$).
Also, write 
\begin{align}
\mr{Mor} (a, b)^\ncong \  \left(\text{resp.},  \ \mr{Mor}\{ a, b\}^\ncong \right)
\end{align}
 for the subset of $\mr{Mor}(a, b)$ (resp., $\mr{Mor}\{ a, b\}$) consisting of non-invertible  morphisms.
  An object  $a$ in $\mcC$ is said to be  {\bf minimal} if  it is not an initial object  and   any monomorphism $b \migiincl a$ in $\mcC$, where $b$ is not an initial object, is necessarily an isomorphism.

Next, we denote by $\mcC^\mr{op}$ the opposite category of $\mcC$.
Also, given  each functor $F : \mcC \migi \mcD$, we shall write $F^\mr{op} : \mcC^\mr{op} \migi \mcD^\mr{op}$ for the functor between the respective opposite categories naturally  induced  by $F$. 
We denote by $\mcC^{\cong}$ the category whose objects are the elements of $\mr{Ob}(\mcC)$ and whose morphisms are the isomorphisms in $\mr{Mor}(\mcC)$.
Given two categories $\mcC$, $\mcD$, we shall write $\mcC \cong \mcD$ (resp., $\mcC \stackrel{\mr{isom}}{\cong} \mcD$) if $\mcC$ is equivalent (resp., isomorphic) to $\mcD$.

Recall that a category $\mcC$ is said to be a groupoid
if the  morphisms in $\mcC$ are all invertible.
By a {\bf groupoid in $\mcC$}, we  mean a subcategory $\mcD$ of $\mcC$ forming a groupoid such that,
if $a, a'$ are objects in $\mcD$  and  $u : a \isom a$ is an isomorphism in  $\mcC$, then $u$ lies  in $\mr{Mor}(\mcD)$.

A category $\mcC$ is said to be a {\bf preorder} if $\mr{Mor}(a, b)$ has at most one element for all objects $a, b \in \mr{Ob}(\mcC)$.
We write $a < b$ if $\mr{Mor}(a, b)$ is nonnempty;
the binary relation ``$<$" in $\mr{Ob}(\mcC)$  is reflexive and transitive.
If in addition ``$<$"
is symmetric, $\mcC$ is said to be an {\bf order}.
Note that ``$<$"  is symmetric precisely when the only isomorphisms of the preorder $\mcC$ are the identity morphisms.
If $(T, <)$ is a partial order, then we write $T^\perp$  for  the  order defined in such a way that $\mr{Ob}(T^\perp) = T$ and $\mr{Mor}(a, b) \neq \emptyset$ precisely when $a < b$.  

Throughout the present paper, we shall fix a universe $\msU$.
Denote by $\mcS et$ the category consisting  of $\msU$-small sets and maps between them.
Denote by $\mcC at$ the category consisting of $\msU$-small categories and  functors between them.
Moreover, denote by $\mcO rd$ the full subcategory of $\mcC at$ consisting of $\msU$-small orders.
To each $\msU$-small set  $T$, we associate the category $\mcD is(T)$, whose objects are the elements in $T$ and whose only morphisms are identity morphisms.
The assignment $T \mapsto \mcD is (T)$ defines a functor $\mcD is : \mcS et \migi \mcO rd$.

\LSP
\subsection{Structure species} \label{SS06}

Let us  recall  the categorical reformulation of Bourbaki's  structure species discussed in ~\cite{Bla}.
Fix a $\msU$-small category $X$.

\SSP
\bde[cf. ~\cite{Bla}, \S\,2, Definition 2.1]  \label{D01}
 Denote by $J : X^{\cong} \migiincl  X$ the natural inclusion.
\begin{itemize}
\item[(i)]
Let us  consider a pair 
\begin{align}
\Sigma := (E, S)
\end{align}
 consisting of two functors $E: X \migi \mcO rd$, $S: X^{\cong}  \migi \mcS et$.
We say that $\Sigma$ is 
a {\bf (covariant) structure species} on $X$
if the composite $\mcD is \circ S$ is a subfunctor of $E \circ J$, meaning that,
for every $a \in  \mr{Ob}(X^{\cong})$, the category $(\mcD is \circ S) (a)$ is a subcategory of $(E \circ J)(a)$ and the inclusion $(\mcD is \circ S) (a) \migiincl (E \circ J)(a)$ is functorial with respect to $a$.
\item[(ii)]
Let $\Sigma := (E, S)$ and $\Sigma' := (E', S')$ be structure species on $X$.
A {\bf morphism of structure species} from $\Sigma$ to $\Sigma'$ is defines as a natural transformation $\phi : S \migi S'$ such that, for every morphism $u : a \migi b$ in $X$ and for every $U \in \mr{Ob}(S (a))$ and  $V \in \mr{Ob}(S(b))$, we have $E' (u)(\phi_a (U)) < \phi_b (V)$  whenever $E (u)(U) < V$.
\end{itemize}
  \ede
\SSP

\begin{rema} \label{R01}
Let $\Sigma := (E, S)$ be a structure species on $X$.
As mentioned in the Remark following ~\cite{Bla}, \S\,2, Definition 2.1,
 $E$ and $S$ respectively  correspond to the {\it echelon} and  {\it structure scheme} of $\Sigma$ in the traditional terminology  of Bourbaki (cf. ~\cite{Bou}, Chap.\,IV, \S\,1).
\end{rema}
\SSP

\begin{exa}[Continuous maps on topological spaces] \label{Ex07}
For each $\msU$-small set $T$, we shall write $P (T)$ for the power set (i.e., the set of subsets) of $T$.
Let $P^+$ (resp., $P^-$) be the functor $\mcS et \migi \mcS et$ (resp., $\mcS et^\mr{op} \migi \mcS et$) 
defined as follows:
\begin{itemize}
\item
 For each object $T$ in $\mcS et$, we set $P^+(T) :=  P(T)$ (resp., $P^- (T) := P(T)$).
\item
 For each morphism $f  : T \migi T'$ in $\mcS et$,   we set $P^+ (f)$ (resp., $P^- (f)$) to be   the map $P^+ (T) \migi P^+(T')$ given by $U \mapsto f(U)$ for every  $U \in P^+ (T)$  (resp., the map $P^- (T') \migi P^-(T)$ given by $U' \mapsto f^{-1}(U')$ for every  $U' \in P (T')$).
\end{itemize}
The maps $P^+ (f)$ and $P^- (f)$ defined for 
 any morphism $f :T \migi T'$ in $\mcS et$
 are  non-decreasing when 
 both $P (T) \left(= P^+ (T) = P^- (T) \right)$ and $P(T') \left(= P^+(T') = P^- (T') \right)$ are equipped with the order structures determined  by ``$\subseteq$".
 Hence, the assignments $T \mapsto (P^+ (T), \subseteq)^\perp$ and $T \mapsto (P^- (T), \subseteq)^\perp$ respectively  induce functors  $\mcP^+ : \mcS et \migi \mcO rd$ and $\mcP^- : \mcS et^{\mr{op}} \migi \mcO rd$.

Now, let $Z$ be  a $\msU$-small subcategory of $\mcS et^\mr{op}$.
We write $E := \mcP^+ \circ P^- |_Z : Z \migi \mcS et \migi \mcO rd$.
For each $T \in \mr{Ob}(Z)$, denote the set of topologies on $T$  by $\mcT op (T)$.
If $f : T \isom T'$ is a bijection of $\msU$-small sets, then $P^+ (P^- (f))$ induces a bijection $\mcT op (f) : \mcT op (T') \isom \mcT op (T)$.
The assignments $T \mapsto \mcT op (T)$ and $f \mapsto \mcT op (f)$ together  yield a functor 
$\mcT op : Z^{\cong} \migi \mcS et$. 
The resulting pair
 \begin{align}
 \Sigma_{\mr{top}} := (E, \mcT op)
 \end{align}
  forms a structure species on $Z$.
\end{exa}
\SSP

We shall denote by
\begin{align}
\mcS p_{X}
\end{align}
the category consisting of structure species on $X$ and morphisms between them.
If there exists an equivalence of categories $X \isom X'$ (where $X'$ is another $\msU$-small category), then we can construct an equivalence of categories $\mcS p_X \isom \mcS p_{X'}$ in an evident manner.

\LSP
\subsection{Constructive functors} \label{SS07}

Next, let us recall an equivalent realization of a structure species, which is  a constructive functor in the sense of   ~\cite{Bla}.
 In ~\cite{Son}, \S\,2, Definition 3, this notion was introduced under the name ``structure species".

\SSP
\bde[cf. ~\cite{Bla}, \S\,3, Definitions 3.1 and 3.2]  \label{D02}
\begin{itemize}
\item[(i)]
Let $F : Y \migi X$ be a functor between $\msU$-small  categories.
We say that $F$ is  a {\bf constructive functor} on $X$ 
(or simply, $F$ is {\bf constructive})
if it satisfies the following conditions:
\begin{itemize}
\item
$F$ is faithful.
\item
For every $a \in \mr{Ob}(Y)$ and for every isomorphism $u$ in $X$ with domain $F (a)$,
there exists uniquely an isomorphism $u_Y$ in $Y$ with domain $a$ satisfying  the equality  $F(u_Y) = u$.
\end{itemize} 
We often refer to the second condition  as the {\it lifting property} on $F$.
\item[(ii)]
Let $F: Y \migi X$ and $F' : Y' \migi X$ be constructive functors on $X$.
A {\bf morphism of constructive functors} from $F$ to $F'$
is defined as a functor $\Phi : Y \migi Y'$ satisfying the equality  $F' \circ \Phi = F$.
\end{itemize}
 \ede
\SSP

We shall denote by
\begin{align}
\mcC on_X
\end{align}
the category consisting of constructive functors on $X$ and morphisms between them. 
One may  verify that $\mcC on_X$ admits finite coproducts and fiber products.
If $f : X' \migi X$ is a functor between $\msU$-small categories, then the assignment $Y \mapsto f^*(Y) := Y \times_X X'$ defines a functor
\begin{align} \label{E100}
f^* :  \mcC on_X \migi \mcC on_{X'}.
\end{align}

Here, we shall construct a constructive functor associated to a structure species.
Let $\Sigma := (E, S)$ be a structure species on $X$.
Denote by $Y_\Sigma$ the category defined as follows:
\begin{itemize}
\item
The objects of $Y_\Sigma$ are pairs $(a, T)$ such that $a \in \mr{Ob}(X)$ and $T \in S (a)$.
\item
The morphisms from $(a, T)$ to $(a', T')$ are morphisms $u : a \migi a'$ in $X$ such that $E (u) (T) < T'$, where ``$<$" is the relation defined on the order $E(a')$.
The composite law for morphisms in $Y_\Sigma$  is defined in a natural manner.
\end{itemize}
Moreover, we set
\begin{align} \label{E12}
F_\Sigma : Y_\Sigma \migi X
\end{align}
to be the functor given by $(a, T) \mapsto a$ and $u \mapsto u$.
This functor  forms a constructive functor on $X$ (cf.  ~\cite{Bla}, \S\,2, Propositions 2.1 and 2.2).
To a morphism of structure species $\phi : \Sigma \migi \Sigma'$, we can associate a morphism of constructive functors $F_\phi : F_\Sigma \migi F_{\Sigma'}$ (but we will omit the details of this construction). 
According to ~\cite{Bla}, \S\,4, Theorem 2, the assignments $\Sigma \mapsto F_\Sigma$ and $\phi \mapsto F_\phi$ together define an equivalence of categories
\begin{align} \label{E14}
\mcS p_X \isom \mcC on_X.
\end{align}

\LSP
\subsection{Property of a constructive functor} \label{SS09}

In this section, we observe one property of  a constructive functor (cf. Proposition \ref{EedP}).
Before  doing so,  let us  recall the connectedness of a category.

\SSP
\bde \label{D03} 
A category $\mcC$  is said to be  {\bf connected} if it is nonempty and  
for any pair of objects $a, b \in \mr{Ob}(\mcC)$, 
there exists a   finite sequence   $(c_1, \cdots, c_n)$ of objects  in $\mcC$ 
 such that $c_1 = a$, $c_n = b$, and  $\mr{Mor} \{ c_j, c_{j+1}\} \neq \emptyset$ for any $j=1, \cdots, n-1$.
\ede
\SSP

\begin{exa} \label{Ex02}
A category $\mcC$ is connected if it has either an initial object or a terminal object.
In particular, a category consisting  of $\msU$-small sets containing either the empty set or  a singleton is connected.
\end{exa}

\SSP

Next, we prove   the  following assertion, which  will be used  in the subsequent discussion.

\SSP
\bpr \label{EedP}
Let $F : Y \migi X$ and  $F' : Y' \migi X$ be  constructive functors on $X$ such that $\mr{Ob}(Y) \neq \emptyset$ and $Y'$ is a connected groupoid.
Also, let $\Phi : F \migi F'$ be a morphism of constructive functors.
Then, the maps  $\mr{Ob}(\Phi) : \mr{Ob}(Y) \migi \mr{Ob}(Y')$ and $\mr{Mor}(\Phi) : \mr{Mor}(Y) \migi \mr{Mor}(Y')$ induced by $\Phi$ are surjective.
\epr
\begin{proof}
By assumption  $\mr{Ob}(Y) \neq \emptyset$,  there exists an object $a_Y$ in $Y$.
Write $a_{Y'} := \Phi (a_Y)$.
Now, let us take an arbitrary  object $b_{Y'}$ in $Y'$.
Since $Y'$ is a connected groupoid, we can find an isomorphism $u_{Y'} : a_{Y'} \isom b_{Y'}$ in $Y'$.
If we write $u := F' (u_{Y'})$, then the lifting property on $F$ implies that 
there exists an isomorphism $u_Y : a_Y \isom b_Y$ in $Y$ (for some $b_Y \in \mr{Ob}(Y)$) with $F(u_Y) = u$.
Since both $u_{Y'}$ and $\Phi (u_Y)$ are  isomorphisms with domain $a_{Y'}$ lifting $u$,
the lifting property on $F$ again implies that $u_{Y'} =\Phi (u_Y)$ and hence $b_{Y'} = \Phi (b_Y)$.
This shows the surjectivity of $\mr{Ob}(\Phi)$.
The surjectivity of $\mr{Mor}(\Phi)$ can be proved by using  a similar argument, so we will finish the proof here.
\end{proof}

\LSP
\subsection{Constructive functor on a groupoid} \label{SS09}

This subsection deals with 
a specific  constructive functor on a groupoid.
From Proposition \ref{P014} described later, 
 this constructive functor 
   can be  thought of as  the universal covering of a topological space or a graph (cf. Example \ref{R04}, Remark \ref{Err0}).

Let  $G$ be  a   $\msU$-small groupoid and $e$ an object in $G$.
 Suppose further  that $G$ is connected, which implies  $\mr{Mor}(a, a') \neq \emptyset$ for any $a, a' \in \mr{Ob}(G)$.
We shall set  $Y_{G, e}$ to be the category defined  as follows:
\begin{itemize}
\item
The objects in $Y_{G, e}$  are  pairs $(a, u)$ consisting of $a \in \mr{Ob}(G)$ and $u \in \mr{Mor}(e, a)$.
\item
The morphisms in $Y_{G, e}$ from $(a, u)$ to $(a', u')$
 are  morphisms   $v : a \migi a'$ satisfying   
$u' = v \circ u$.
The composition law for  morphisms in $Y_{G, e}$ is defined in an evident manner.
\end{itemize}
One may verify that $Y_{G, e}$ is a connected groupoid and,   for any $a, a' \in \mr{Ob}(Y_{G, e})$, the set  $\mr{Mor} (a, a')$ has exactly one element.
The assignments $(a, u) \mapsto a$ and $v \mapsto v$ together  define  a functor 
\begin{align} \label{E012}
F_{G, e} : Y_{G, e} \migi G,
\end{align}
 which is verified to form a constructive functor on $G$.

\SSP
\begin{exa} \label{R04}
Let us consider the case where $\mr{Ob}(G)$ has exactly one element, which we denote by  $\circledast$.
Denote by $\mr{Aut}(\circledast)$ the automorphism  group of  $\circledast$.
Then, the category $Y_{G, \circledast}$ may be  regarded as 
 the Cayley graph $\mr{Cay} (\mr{Aut}(\circledast))$ of  the  group $\mr{Aut}(\circledast)$ by taking the connection set  as $\mr{Aut}(\circledast)$ itself, where   the  vertices and arcs in $\mr{Cay} (\mr{Aut}(\circledast))$ are  respectively associated  with the objects and morphisms in $Y_{G, \circledast}$.

\end{exa}
\SSP

We shall prove  the following proposition concerning   the constructive functor $F_{G, e}$.

\SSP
\bpr \label{P014} 
Let $F: Y \migi G$ be a constructive functor on $G$. 
\begin{itemize}
\item[(i)]
Denote by $F^{-1}(e)$ the preimage of  the subset $\{ e\} \left( \subseteq \mr{Ob}(G)\right)$ via the map $\mr{Ob}(Y) \migi \mr{Ob}(G)$ induced by $F$.
Then, the  map of sets
\begin{align}
\mr{ev}_{F, e} : \mr{Mor} (F_{G, e}, F) \isom F^{-1}(e)
\end{align}
obtained  by assigning $\Phi \mapsto \Phi ( (e, \mr{id}_e))$ is bijective.
\item[(ii)]
Suppose that $Y$ is a connected groupoid.
Then, any morphism  of constructive functors $\Phi : F \migi F_{G, e}$ is an isomorphism.
In particular, any endomorphism of $F_{G, e}$ is an isomorphism.
\item[(iii)]
Suppose that $\mr{Ob}(Y) \neq \emptyset$.
Then, any monomorphism of constructive functors $\Phi : F \migiincl F_{G, e}$ is an isomorphism.
\end{itemize}
 \epr
\begin{proof}
First,  we prove the surjectivity of $\mr{ev}_{F, e}$ in assertion  (i).
Let $e_Y$ be an element of  $F^{-1} (e)$.
In what follows, we shall construct a morphism of constructive functors  $F_{G, e} \migi F$ that  is mapped to $e_Y$ via $\mr{ev}_{F, e}$.
Let us take an arbitrary object $(a, u)$ of $Y_{G, e}$.
Since $F$ is constructive, there exists a unique pair $(a_Y, u_Y)$ consisting of an object $a_Y$ in $Y$ and a morphism $u_Y : e_Y \migi a_Y$ with $F (u_Y) = u$.
Next, let $v: (a, u) \migi (a', u')$ be a morphism in $Y_{G, e}$.
Denote by $v_Y$ the unique lifting of $v$ with domain $a_Y$.
Since $F$ is constructive,  the morphism $u' \left(= v \circ u\right)$ lifts uniquely to a morphism in $Y$ with domain $e_Y$.
This implies  $a'_Y = v_Y \circ a_Y$.
Moreover, the lifting property on $F$ again implies that $(\mr{id}_{a})_Y = \mr{id}_{a_Y}$ and  $(v' \circ v)_Y = v'_Y \circ v_Y$ for any morphism
 $v' :(a', u') \migi (a'', u'')$.
 Thus, the assignments $(a, u) \mapsto a_Y$ and $v \mapsto v_Y$ define a morphism
 $\Phi_{e_Y} : F_{G, e} \migi F$, and this morphism satisfies $\mr{ev}_{F, e}(\Phi_{e_Y}) = e_Y$ by construction.
 This implies the surjectivity of $\mr{ev}_{F, e}$.
 
 Next, we prove the injectivity of $\mr{ev}_{F, e}$.
 Let $\Phi$ be a morphism $F_{G, e} \migi F$ and write $e_Y := \mr{ev}_{F, e} (\Phi)$.
 The problem reduces  to one of  proving  the equality  $\Phi = \Phi_{e_Y}$.
 To this end, let us take an arbitrary element $(a, u)$ of $Y_{G, e}$.
Since $u$ defines a morphism $(e, \mr{id}_e) \migi (a, u)$ in $Y_{G, e}$,
the image $\Phi (u)$ defines a morphism $e_Y \left(= \Phi ((e, \mr{id}_e)) \right) \migi \Phi ((a, u))$.
This morphism is a unique lifting of $u \in \mr{Mor}(G)$ with domain $e_Y$, so we have $\Phi ((a, u)) = \Phi_{e_Y}((a, u))$.
Thus, we obtain  the equality  $\Phi = \Phi_{e_Y}$, as desired.
This completes the proof of   the bijectivity  of $\mr{ev}_{F, e}$.

Now, we prove  assertion (ii).
From Proposition \ref{EedP} and  the fact that $Y_{G, e}$ is a connected groupoid, it suffices to prove the injectivity of $\mr{Ob}(\Phi)$ and of  $\mr{Mor}(\Phi)$.
Suppose that there exists two objects  $c_Y$, $c'_Y$  in $Y$ with $\Phi (c_Y) = \Phi (c'_Y) \left(=:(c, v) \right)$.
Since $Y$ is a connected groupoid, there exists  an isomorphism $v_Y : c_Y \isom c'_Y$ in $Y$.
The set $\mr{Mor}((c, v), (c, v))$ coincides with $\{\mr{id}_v \}$, so the equality $\Phi (v_Y) = \mr{id}_v$ holds.
The equality  $F = F_{G, e} \circ \Phi$ implies $F (v_Y ) = \mr{id}_Y$.
But, by the lifting property on $F$, $v_Y$ must be equal to $\mr{id}_{c_Y}$, which implies $c_Y = c'_Y$.
This completes the proof of the injectivity of $\mr{Ob}(\Phi)$.
A similar argument can be used to proved the
 injectivity of $\mr{Mor}(\Phi)$.

Finally, we prove assertion (iii).
Let $\Phi : F \migiincl F_{G, e}$ be a monomorphism. 
Since $Y \neq \emptyset$, there exists a connected groupoid  $G_0$ in  $Y$.
The restriction $F |_{G_0} : G_0 \migi G$ of $F$ to $G_0$ forms  a constructive functor on $G$.
Then, it follows from assertion (ii) that
the  morphism   $\Phi |_{G_0} : F |_{G_0} \migi F_{G, e}$ obtained by restricting  $\Phi$ is an isomorphism.
Here, suppose that there exists an object $a \in \mr{Ob}(Y) \setminus \mr{Ob}(G_0)$.
Denote by $G_1$ the groupoid in $Y$ containing $a$ (hence $\mr{Ob}(G_0) \cap \mr{Ob}(G_1) = \emptyset$).
For the same reason as above, the restriction $\Phi |_{G_1} : F |_{G_1} \migi  F_{G, e}$ of $\Phi$ to $G_1$ forms an isomorphism of constructive functors.
Thus, we obtain  two distinct morphisms $(\Phi |_{G_0})^{-1}$, $(\Phi |_{G_1})^{-1}$ belonging to  $\mr{Mor}(F_{G, e}, F)$  that coincide with each other after composing  with $\Phi$.
This  contradicts the assumption that $\Phi$ is a monomorphism.
This implies  $\mr{Ob}(G_0) = \mr{Ob}(Y)$.
Next, suppose that there exists a morphism $u \in \mr{Mor}(Y)\setminus \mr{Mor}(G_0)$.
Then, the two morphisms $\mr{id}_Y$ and $(\Phi |_{G_0})^{-1} \circ \Phi$ are distinct but   identical  to each other when composed  with $\Phi$.
This is a contradiction, so we have $\mr{Mor}(Y) = \mr{Mor}(G_0)$.
 This  implies that $Y = G_0$ and $\Phi = \Phi |_{G_0}$.
Thus,  $\Phi$ turns out to be an isomorphism.
\end{proof}
\SSP

Let $e'$ be an object in $G$
  and $w : e \migi e'$ a morphism in $G$.
Then,  
we obtain a morphism of constructive functors
\begin{align} \label{E1013}
\Phi_w := \mr{ev}_{F_{G, e'}}^{-1}((e, w^{-1})) : F_{G, e} \migi  F_{G, e'}.
\end{align}
To be precise, this morphism is obtained  by assigning
  $(a, u) \mapsto (a, u \circ w^{-1})$ (for each $(a, u) \in \mr{Ob} (Y_{G, e})$) and  $v \mapsto v$ (for each $v \in \mr{Mor}(Y_{G, e})$).
 From Proposition \ref{P014}, (ii), $\Phi_w$ turns out to be  an isomorphism.
The existence of $\Phi_w$ implies that  the isomorphism class of  the constructive functor  $F_{G, e}$ does not depend on the choice of $e$.
By considering the case of  $e = e'$ and applying again  
Proposition \ref{P014}, (ii),
we see that
the assignment $w \mapsto \Phi_w$ defines an isomorphism  of groups
\begin{align} \label{E1234}
\phi_e := \mr{ev}_{F_{G, e}, e}^{-1} : \left(F^{-1}_{G, e} (e) =  \right) \mr{Aut} (e) \migi \mr{Aut}_{\mcC on_G} (F_{G, e}),
\end{align}
where 
$\mr{Aut}_{\mcC on_G} (F_{G, e})$ denotes the automorphism group of $F_{G, e}$ in $\mcC on_{G}$.

\SSP
\begin{rema} \label{Err0}
One may verify   that $\mcC on_G$ forms a Galois category  (cf. ~\cite{Gro1}, Expos\'{e} V,  Th\'{e}or\`{e}m 4.1 and D\'{e}finition 5.1)  by setting   $\mr{ev}_{F, e}$ as the fiber functor.
The above discussion shows that the resulting  fundamental group $\pi_1 (\mcC on_G; \mr{ev}_{F, e})$ of $\mcC on_G$ with base point $\mr{ev}_{F, e}$  is isomorphic to the group $\mr{Aut}(e)$.
Hence, the equivalence class of the groupoid $G$ can be reconstructed from the fundamental group $\pi_1 (\mcC on_G; \mr{ev}_{F, e})$.
This fact is reminiscent of ``Grothendieck  conjecture"-type theorems in anabelian geometry  (cf. Remark \ref{R450}).
\end{rema}

\SSP
\begin{rema} \label{Err1}
The cultural model, as well as  the paradigm for knowledge and practice,  in the modern Western world envisioned by the structuralism of L\'{e}vi-Strauss and others is characterized by an organizational structure of an {\it arborescence} system that looks for the origin of ``things" and for the culmination or conclusion of those ``things".
In this model, a small idea, like a seed, 
takes root and grows into a tree with a sturdy trunk supporting numerous branches, all linked to and traceable back to the original idea. 
\footnote{Due to criticism of how an object can be explained  only with such a picture,
 a perspective that focuses (not only on the origin and conclusion but)  on the process of its creation and the possibility of its change arose later. It led to the establishment of the position known today as {\it post-structuralism}.  
 In this context,
 G. Deleuze and F. Guattari used the term ``{\it rhizome}" to describe a process of existence and growth that does not come from a single central point of origin (cf. ~\cite{DeGu}). }
The results of Proposition \ref{P014} together with the  viewpoint of the previous  Remark  suggest that, in accordance with the picture of  L\'{e}vi-Strauss' structuralism, 
$F_{G, e}$ plays the role of a seed or a trunk  in the system of structure species on $G$.
\end{rema}

\vspace{10mm}
\section{Proof of the main theorem}  \label{S44}
\SSP

This section is devoted to proving Theorem \ref{E223}.
 The nontrivial portion of that theorem is  the implication (b) $\Rightarrow$ (a), which asserts 
that
the equivalence class  of a category  $X$ may be characterized uniquely, up to a certain indeterminacy,  from  the categorical structure of  $\mcS p_X$, or equivalently, of $\mcC on_X$ (cf. (\ref{E14})).
In the following discussion,
we will often speak of various
things concerning
$\mcC on_X$ as being ``{\it reconstructed (or characterized) category-theoretically}".
By this, we mean that
they are preserved by an  arbitrary equivalence of categories $\mcC on_X \isom \mcC on_{X'}$ (where $X'$ is another $\msU$-small category).
For instance, the set of monomorphisms in $\mcC on_X$ may be {\it characterized category-theoretically} as the morphisms $\Phi : F \migi F'$ such that,
for any $F''$, the map of sets $\mr{Map}(F'', F) \migi \mr{Map}(F'', F')$ obtained  by composing with $\Phi$ is injective.
To simplify our notation, however, we will omit explicit mention of this equivalence $\mcC on_X \isom \mcC on_{X'}$, of $X'$, and of the various ``primed" objects and morphisms corresponding to the original objects and morphisms, respectively, in $\mcC on_X$.
 
Our tactic for  completing the proof of the implication (b) $\Rightarrow$ (a)  (i.e., recognizing the structure of $X$)
is, as in ~\cite{Moc1}, ~\cite{Moc2}, and ~\cite{Wak}, to reconstruct step-by-step various partial information about  $X$ from the categorical structure of $\mcC on_X$.

\LSP
\subsection{Reconstruction of the objects and their automorphisms} \label{SS57}
\SSP

The first step of the proof is to 
reconstruct the set of isomorphism classes of objects in $\mcC on_X$ and their automorphism groups.
Let us fix a $\msU$-small category $X$.
Also, let us fix a skeleton $\overline{X}$ of $X$ (i.e., a full subcategory $\overline{X}$ of $X$ such that the inclusion $\overline{X} \migiincl X$ is an equivalence of categories  and  no two distinct objects in $\overline{X}$ are isomorphic).

If $G$ and  $G'$ are distinct connected groupoids in $X$, then we see that $\mr{Ob}(G) \cap \mr{Ob}(G') = \emptyset$.
Moreover, for an object $e$  of $X$,  there exists a unique connected groupoid $G_e$ in $X$ containing $e$.
Thus, the assignment $e \mapsto G_e$ defines  a  bijective correspondence  between $\mr{Ob}(\overline{X})$ and the set of connected groupoids in $X$.
Also,  
we have a decomposition
\begin{align}
\mr{Ob} (X) = \coprod_{G} \mr{Ob}(G),
\end{align}
where the disjoint union in the right-hand side runs over the set of connected groupoids in $X$.

Now, let $e$ be an object in $\overline{X}$.
 Denote by 
$G$ the connected groupoid in $X$ containing $e$ (hence $\mr{Ob}(\overline{X}) \cap \mr{Ob}(G) = \{e \}$) and 
 denote the natural inclusion by $\iota_G : G \migiincl X$.
The assignment $F \mapsto  \iota_{G*}(F) := \iota_G \circ F$ defines a fully faithful functor
\begin{align} \label{E1005}
\iota_{G*} : \mcC on_G \migi \mcC on_X.
\end{align}
The following proposition can be proved immediately from the definitions of a constructive functor and $\iota_{G*}$ (so we will omit their proofs.)

\SSP
\bpr \label{PP01}
\begin{itemize}
\item[(i)]
This functor is left adjoint to the functor $\iota^*_G$ (cf. (\ref{E100})), which 
 means that
for  $F \in \mr{Ob}(\mcC on_X)$ and  $F' \in \mr{Ob}(\mcC on_G)$, there exists a functorial bijection
\begin{align} \label{E1004}
\mr{Mor}(\iota_{G*} (F'), F) \isom \mr{Mor} (F', \iota_G^* (F)).
\end{align}
\item[(ii)]
Let $F$ be a constructive functor on $X$.
Suppose that there exist a constructive functor $F'$ on $G$ and  a morphism $F \migi \iota_{G*}(F')$.
Then, there exists a constructive functor $F''$ on $G$ with $\iota_{G*}(F'') \cong F$.
\end{itemize}
\epr
\SSP

The following assertion is a direct consequence of assertion (ii) above.

\SSP
\bco \label{CC01}
Let $F$ be a constructive functor on $G$.
Then, $F$ is minimal in $\mcC on_G$ if and only if $\iota_{G*}(F)$ is minimal in  $\mcC on_X$.
\eco
\SSP

We shall write
\begin{align} \label{E1006}
F^+_{G, e} :=  \iota_{G*} (F_{G, e}) : Y_{G, e}  \migi  X.
\end{align}
When there is no fear of confusion, we will write $F^+_{G} := F^+_{G, e}$ for simplicity.
 (This abbreviation of the notation can be justified because of  the existence of  (\ref{E1013})  and  the fully faithfulness of $\iota_{G*}$.
In fact, these facts show that  the isomorphism  class of $F^+_{G, e}$ depends only on $G$; i.e.,  it does not depend on the choice of the skeleton $\overline{X}$.)

Denote by 
\begin{align}
{^\dagger}\mr{Ob}(\overline{X}) \left(\subseteq \overline{\mr{Ob}} (\mcC on_X) \right)
\end{align}
the set of isomorphism classes of objects in $\mcC on_X$ of the form $F^+_G$ for  some  connected groupoid $G$ in $X$.
If $G$ and $G'$ are connected groupoids in $X$, then $F^+_G \cong F^+_{G'}$ precisely when $G = G'$.
This implies that the assignment
 from each  element $e' \in \mr{Ob}(\overline{X})$ to the constructive functor $F_{G'}^+$, where $G'$ denotes the connected groupoid in $X$ containing $e'$, defines a bijection of sets
\begin{align}
\xi_{\overline{X}} : \mr{Ob}(\overline{X}) \isom {^\dagger}\mr{Ob}(\overline{X}). 
\end{align}

\SSP
\bpr \label{L01} 
Let $F$ be a constructive functor on $X$.
Then, $F$ is minimal in $\mcC on_X$ if and only if $F$ is isomorphic to  $F^+_G$ for some connected groupoid $G$ in $X$.
In particular, the subset 
${^\dagger}\mr{Ob}(\overline{X})$ of $\overline{\mr{Ob}}(\mcC on_X)$
can be reconstructed category-theoretically from the data $\mcC on_X$ (i.e., of a category). 
 \epr
\begin{proof}
The ``if" part of the required equivalence   follows immediately from Proposition \ref{P014}, (iii), and Corollary \ref{CC01}.
Hence, in what follows, we shall prove
 the ``only if" part.
 Let $F : Y \migi X$ be a minimal object in $\mcC on_X$.
 Since this object is nonempty, there exists an object $a$ in $Y$.
 Then, by the bijection $\mr{ev}_{F, F (a)}$ in Proposition \ref{P014}, (i), together with the adjunction relation (\ref{E1004}),
 we can find a morphism of constructive functors $\Phi : F_{G}^+ \migi F$.
 It follows from Proposition \ref{EedP}
 that the maps $\mr{Ob}(\Phi) : \mr{Ob}(Y_{G, e}) \migi \mr{Ob}(Y)$, $\mr{Mor}(\Phi) : \mr{Mor}(Y_{G, e}) \migi \mr{Ob}(Y)$ induced by $\Phi$ are surjective.
 Next,  suppose that we are given  two objects $(a, u)$, $(a, u')$ in $Y_{G, e}$ with $\Phi ((a, u)) = \Phi ((a, u'))$.
 After possibly composing $\Phi$ with a suitable automorphism of $F_G^+$,
 we may assume that $a = e$.
 Denote by $v$ the unique endomorphism of $e$ in $G$  satisfying  $\{ v \} = \mr{Mor}((e, u), (e, u'))$.
 Then, from the bijectivities of $\mr{ev}_{F, e}$ and $\mr{ev}_{F_{G, e}}$,
 it can immediately be seen that $v = \mr{id}_e$, i.e., $(e, u) = (e, u')$.
  This implies the injectivity of  $\mr{Ob}(\Phi)$.
 Moreover, since $\mr{Mor}((a, u), (b, v))$ is a singleton for any objects $(a, u)$, $(b, v)$ in $Y_{G, e}$, the injectivity of $\mr{Ob}(\Phi)$ implies that of $\mr{Mor}(\Phi)$.
 Consequently,
  $\Phi$ turns out to be an isomorphism.
 This completes the proof of the ``only if" part.
 \end{proof}
\SSP

Let $e$ and $G$ be as above.
Then,  for each $w \in \mr{Aut}(e)$, the automorphism $\Phi_w$ (cf. (\ref{E1013})) defines, via $\iota_{G*}$, an automorphism $\Phi^+_u$ of $F^+_G$.
By  the bijectivity of $\phi_e$ (cf. (\ref{E1234})) and   the fully faithfulness of $\iota_{G*}$,
 the assignment $w \mapsto \Phi_w^+$ 
 defines an isomorphism of groups
\begin{align} \label{E1235}
\xi^{\cong}_e : \mr{Aut} (e) \isom {^\dagger}\mr{Aut} (e),
\end{align}
where  we set $ {^\dagger}\mr{Aut} (e) := \mr{Aut}_{\mcC on_X}(F^+_G)$.
The following assertion can be verified immediately.

\SSP
\bpr \label{L090}
The  subset ${^\dagger}\mr{Aut} (e)$ of $\mr{Mor}(\mcC on_X)$ together with its group structure can be reconstructed category-theoretically from the data $(\mcC on_X, F^+_{G})$ (i.e., of a category and a minimal object in this category).
\epr

\LSP
\subsection{Reconstruction of the non-invertible morphisms} \label{SS58}
\SSP

Let us take two distinct objects $e$, $e'$ in $\overline{X}$.
Denote by $G$ and $G'$ the connected groupoids in $X$ containing $e$ and $e'$ respectively.
Given a  non-invertible morphism $v \in \mr{Mor}(e, e')$ (i.e., $v \in \mr{Mor}(e, e')^\ncong$), 
we shall denote by
\begin{align}
Y_{G, G', v},  \ \text{or simply} \ Y_v, 
\end{align}
the category determined  as follows:
\begin{itemize}
\item
$\mr{Ob}(Y_{G, G', v} ) = \mr{Ob}(Y_{G, e}) \sqcup \mr{Ob}(Y_{G', e'})$.
\item
$\mr{Mor}(Y_{G, G', v} ) = \mr{Mor}(Y_{G, e}) \sqcup \mr{Mor}(Y_{G', e'}) \sqcup \coprod_{\substack{(a, u) \in \mr{Ob}(Y_{G, e}), \\ (a', u') \in \mr{Ob}(Y_{G', e'})}} \mr{Mor}((a, u), (a', u'))$, where
$\mr{Mor} ((a, u), (a', u')) := \{ u' \circ v \circ u^{-1}\}$ for each
$(a, u) \in \mr{Ob}(Y_{G, e})$, $(a', u') \in \mr{Ob}(Y_{G', e'})$.
The composition law for morphisms in $Y_{G, G', v}$ is defined in a natural manner.
\end{itemize}
Moreover,  denote by 
\begin{align}
F^+_{G, G', v} : Y_{G, G', v} \migi X
\end{align}
the functor given by $(a, u) \mapsto a$ (for each $(a, u) \in \mr{Ob}(Y_{G, G', v})$) and $w \mapsto w$ for each $w \in \mr{Mor}(Y_{G, G', v})$.
Then,  $F^+_{G, G', v}$ forms a constructive functor on $X$, and 
the natural inclusions $F^+_{G, e} \migiincl  F^+_{G, G', v}$ and $F^+_{G', e'} \migiincl  F^+_{G, G', v}$ yields 
   a  non-invertible monomorphism 
\begin{align}
\Upsilon_{G, G', v} : F^+_{G, e} \sqcup F^+_{G', e'}\migiincl F^+_{G, G', v}.
\end{align}
When there is no fear of confusion, we will write $F_v^+ := F_{G, G', v}^+$ and $\Upsilon_v := \Upsilon_{G, G', v}$.

Next, denote by 
\begin{align}
\mcC_{e, e'}
\end{align}
 the category defined as follows:
\begin{itemize}
\item
The objects in $\mcC_{e, e'}$ are  pairs 
$(F, \Upsilon)$ consisting of a constructive functor $F: Y \migi X$  on $X$ and a non-invertible 
monomorphism $\Upsilon : F^+_{G, e} \sqcup F^+_{G', e'} \migiincl F$ in $\mcC on_X$.
\item
The morphisms from $(F, \Upsilon)$ to $(F', \Upsilon')$ are  morphisms of constructive functors  $\Psi : F \migi F'$ satisfying  $\Upsilon' = \Psi \circ \Upsilon$.
\end{itemize}
For 
 each non-invertible  morphism $v : e \migi e'$ in $\overline{X}$,
the  pair  $(F^+_{G, G', v}, \Upsilon_{G, G', v})$ introduced  above specifies  an object of $\mcC_{e, e'}$.
The assignment $(F, \Upsilon) \mapsto F$ defines a functor
\begin{align} \label{E477}
\mcC_{e, e'} \migi \mcC on_X.
\end{align}
By taking account of Proposition  \ref{L01}, we see that 
 the category $\mcC_{e, e'}$ together   with the  functor (\ref{E477}) can be reconstructed category-theoretically from  the data $(\mcC on_X, F^+_{G, e}, F^+_{G', e'})$ (i.e., of a category and two  distinct minimal objects in this category).
 We identify $\mcC_{e, e'}$  with $\mcC_{e', e}$ via the equivalence of categories $\mcC_{e, e'} \isom \mcC_{e', e}$ obtained  by switching the factors $F^+_{G, e} \sqcup F^+_{G', e'} \isom F^+_{G', e'} \sqcup F^+_{G, e}$.
 Under this identification,  $F_{G', G, v'}^+$ and $\Upsilon_{G', G, v'}$ (where $v' \in \mr{Mor}(e', e)^\ncong$) may be regarded as  elements  of $\mr{Ob}(\mcC_{e, e'})$ and $\mr{Mor}(\mcC_{e, e'})$ respectively.

Moreover, we  denote by 
\begin{align}
{^\dagger}\mr{Mor} \{e, e'\}^\ncong
\left(\subseteq \mr{Mor} (\mcC on_X) \right)
\end{align}
the set of morphisms in $\mcC on_X$  of the form $\Upsilon_v$ for some $v \in \mr{Mor}\{e, e' \}^\ncong$.
Since  $v \neq v'$ implies  $(F^+_{v}, \Upsilon_{v}) \ncong (F^+_{v'}, \Upsilon_{v'})$,
the assignment $v \mapsto \Upsilon_v$ defines a bijection of sets
\begin{align} \label{E9898}
\xi_{e, e'}^\ncong :  \mr{Mor}\{e, e' \}^\ncong \isom {^\dagger}\mr{Mor} \{e, e'\}^\ncong.
\end{align}

\SSP
\bpr \label{L06}
Let $(F, \Upsilon)$ be an object in $\mcC_{e, e'}$.
Then, 
the following two conditions are equivalent to each other:
\begin{itemize}
\item[(a)]
$(F, \Upsilon)$ is minimal in $\mcC_{e, e'}$ and
there is no triple of morphisms  
\begin{align} \label{E231}
(\Upsilon_0 : F_{G, e}^+ \migi F_0, \Upsilon'_0 : F_{G', e'}^+ \migi F'_0, \Phi : F_0 \sqcup F'_0 \migi F)
\end{align}
 in $\mcC on_X$ such that $\Phi$ is an isomorphism and satisfies the equality $\Upsilon = \Phi \circ (\Upsilon_0 \sqcup \Upsilon'_0)$.
  \item[(b)]
 $(F, \Upsilon)$ is  isomorphic to $(F^+_{v}, \Upsilon_{v})$
 for some  $v \in \mr{Mor} \{e, e'\}^\ncong$.
 \end{itemize}
In particular,  
the subset ${^\dagger}\mr{Mor} \{e, e'\}^\ncong$ of $\mr{Mor}(\mcC on_X)$
  can be reconstructed category-theoretically from the data $(\mcC on_X, F_{G,e}^+, F_{G', e'}^+)$ (i.e., of a category and two distinct  minimal  objects in this category).
  \epr
\begin{proof}
Since the implication (b) $\Rightarrow$ (a) can be immediately verified from the definition of $(F^+_{v}, \Upsilon_{v})$,
we only consider the inverse direction.
Let $(F : Y \migi X, \Upsilon)$ be  an object in $\mcC_{e, e'}$ satisfying the conditions in (a).
Denote by $\overline{Y}_{G, e}$ and $\overline{Y}_{G', e'}$ the subcategories in $Y$ defined as the images via $\Upsilon$ of $Y_{G, e}$ and $Y_{G', e'}$ respectively.
By the lifting property on $F$, $\overline{Y}_{G, e}$ and $\overline{Y}_{G', e'}$ respectively  specify connected groupoids in $Y$.
Now, suppose that there exists 
an object $a_0$ in $\mr{Ob}(Y) \setminus (\mr{Ob}(\overline{Y}_{G, e}) \sqcup \mr{Ob}(\overline{Y}_{G', e'}))$.
If $G_0$ denotes
the connected groupoid in $Y$ containing $a_0$,
then  the minimality of $(F, \Upsilon)$ implies that  $Y$ must be equal to the disjoint union $\overline{Y}_{G, e} \sqcup \overline{Y}_{G', e'} \sqcup G_0$.
But this is a contradiction because of the second condition in (a).
Hence, we have $\mr{Ob}(Y) =  \mr{Ob}(\overline{Y}_{G, e}) \sqcup \mr{Ob}(\overline{Y}_{G', e'})$.
Since $\Upsilon$ is a monomorphism, the functors $Y_{G, e} \migi \overline{Y}_{G, e}$, $Y_{G', e'} \migi \overline{Y}_{G', e'}$ obtained by restricting $\Upsilon$  are an isomorphism (cf. the proof of Proposition \ref{L01});
by using these isomorphisms, we consider $Y_{G, e} \sqcup Y_{G', e'}$ to be  a subcategory of $Y$.
From the second condition in (a),
there exists objects $a \in \mr{Ob}(\overline{Y}_{G, e})$, $a' \in \mr{Ob}(\overline{Y}_{G', e'})$ such that $\mr{Mor}\{ a, a' \}\neq \emptyset$.
Hence,  the subcategory $Y_{G, e} \sqcup Y_{G', e'}$ of 
$Y$ extends to a subcategory  of the form $Y_{G, G',  v}$ for some $v \in \mr{Mor}\{e, e' \}$.
It follows from the minimality of $(F, \Upsilon)$ that 
the inclusion $Y_{G, G', v} \migiincl Y$ must be an isomorphism.
This completes the proof of the implication (a) $\Rightarrow$ (b).
\end{proof}

\LSP
\subsection{Reconstruction of the composition law} \label{SS59}
\SSP

Let $e$, $e'$, $G$,  and $G'$ be as above.
We shall set 
\begin{align}
\mr{Com}(e, e')^{\cong} &:=\left\{ (u, u', u'') \in  \mr{Aut}(e) \times \mr{Mor} (e, e')^\ncong \times  \mr{Mor} (e, e')^\ncong \, \Big| \, u'' = u' \circ u  \right\},\\
\mr{Com}(e, e')^{\cong}_\mr{op} &:=\left\{ (u, u', u'') \in   \mr{Aut}(e) \times  \mr{Mor} (e', e)^\ncong \times   \mr{Mor} (e', e)^\ncong \, \Big| \, u'' = u \circ u'  \right\}.
 \notag
\end{align}
We will prove  the following assertion.

\SSP
\bpr \label{L34}
Let  us choose $u \in \mr{Aut}(e)$ and  $v_1, v_2  \in \mr{Mor} \{ e, e' \}^\ncong$.
Then, the following two conditions are equivalent to each other:
\begin{itemize}
\item[(a)]
Either $v_2 = v_1 \circ u$ or $v_2 = u \circ v_1$ holds.
\item[(b)]
Either $\Upsilon_{v_2} = \Upsilon_{v_1} \circ (\mr{id}_{F_G^+} \sqcup \Phi_{u^{-1}})$ or 
$\Upsilon_{v_2} = \Upsilon_{v_1} \circ (\mr{id}_{F_G^+} \sqcup \Phi_{u})$ holds.
\end{itemize}
In particular 
the image of $\mr{Com}(e, e')^{\cong} \sqcup \mr{Com}(e, e')^{\cong}_\mr{op}$ via the bijection
\begin{align}
\xi^{\cong}_{e, e, e'} := \xi^{\cong}_e \times \xi_{e, e'}^\ncong \times \xi_{e, e'}^\ncong : & \ \mr{Aut}(e) \times \mr{Mor} \{ e,e' \}^\ncong \times \mr{Mor} \{ e,e' \}^\ncong \\
& \  \isom {^\dagger}\mr{Aut}(e) \times {^\dagger}\mr{Mor} \{ e,e' \}^\ncong \times {^\dagger}\mr{Mor} \{ e,e' \}^\ncong \notag
\end{align}
can be reconstructed category-theoretically from the data $(\mcC on_X, F_{G, e}^+, F_{G', e'}^+)$, i.e., of a category and two distinct minimal objects in this category
 (cf. Proposition \ref{L090} for the category-theoretic reconstruction of ${^\dagger}\mr{Aut}(-)$).
\epr
\begin{proof}
The assertion follows from the various definitions involved.
\end{proof}
\SSP

 Let $G, G', G''$ be connected groupoids in $X$.
 Denote by  $e$, $e'$, and $e''$ the objects of $\overline{X}$ belonging 
 to $G$, $G'$, and $G''$, respectively.
Moreover,   set 
\begin{align}
\mr{Com} (e, e', e'')^{\ncong} &:=\left\{ (u, u', u'') \in  \mr{Mor} (e, e')^\ncong   \times \mr{Mor} (e', e'')^\ncong \times  \mr{Mor} (e, e'')^\ncong \, \Big| \, u'' = u' \circ u  \right\}, \\
\mr{Com} (e, e', e'')^{\ncong}_{\mr{op}} &:=\left\{ (u, u', u'') \in  \mr{Mor} (e', e)^\ncong   \times \mr{Mor} (e'', e')^\ncong \times  \mr{Mor}(e'', e)^\ncong \, \Big| \, u'' = u \circ u'  \right\}.  \notag
\end{align}
We will prove  the following assertion.

\SSP
\bpr \label{L07}
Let us choose $v \in \mr{Mor} \{ e, e'\}^\ncong$,
$v' \in \mr{Mor}\{e', e''\}^\ncong$, and 
$v'' \in \mr{Mor}\{ e, e'' \}^\ncong$.
Then, the following two conditions are equivalent to each other
\begin{itemize}
\item[(a)]
Either $v \circ v'$ or $v' \circ v$ can be defined and one of the equalities $v'' = v \circ v'$, $v'' = v' \circ v$ holds.
\item[(b)]
The colimit $F^+_{\mr{limit}}$  of the diagram
\begin{align} \label{E445}
\vcenter{\xymatrix@C=46pt@R=36pt{
F_{G'}^+ \ar[r]^-{\mr{inclusion}} \ar[d]_-{\mr{inclusion}} & F_{G', G'', v'}^+ \\
 F_{G, G', v}^+ &
}}
\end{align}
exists in $\mcC on_X$, and   there exists a monomorphism 
$F^+_{G, G'', v''} \migiincl  F^+_{\mr{limit}}$ in $\mcC on_X$
via which the morphism $\Upsilon_{v''} : F^+_{G, e} \sqcup F^+_{G'', e''} \migiincl F^+_{v''}$ is compatible 
with the morphism $F^+_{G, e} \sqcup F^+_{G'', e''} \migiincl F^+_\mr{limit}$ arising naturally from $\Upsilon_{v}$ and $\Upsilon_{v'}$.
\end{itemize}
In particular,
the image of  $\mr{Com} (e, e', e'')^{\ncong} \sqcup \mr{Com} (e, e', e'')^{\ncong}_{\mr{op}}$ via the bijection
\begin{align}
\xi_{e, e', e''}^{\ncong} := \xi_{e, e'}^{\ncong} \times \xi_{e', e''}^{\ncong} \times \xi_{e, e''}^{\ncong} : & \ \mr{Mor}\{ e, e' \}^\ncong  \times  \mr{Mor} \{ e', e'' \}^\ncong  \times \mr{Mor} \{ e, e'' \}^\ncong  \\
& \  \isom {^\dagger}\mr{Mor}\{ e, e' \}^\ncong  \times  {^\dagger}\mr{Mor} \{ e', e'' \}^\ncong  \times {^\dagger}\mr{Mor} \{ e, e'' \}^\ncong \notag 
\end{align}
can be reconstructed category-theoretically from the data $(\mcC on_X, F^+_{G, e}, F^+_{G', e'}, F^+_{G'', e''})$, i.e., a category and three distinct minimal objects in this category
  (cf. Proposition \ref{L06} for the category-theoretic reconstruction of ${^\dagger}\mr{Mor}\{ -, - \}^\ncong$).
 \epr
\begin{proof}
The assertion follows from the various definitions involved.
\end{proof}
\SSP

\begin{rema} \label{R340}
Let $a \in \mr{Mor}(e,  e')$ and $b \in \mr{Mor}\{ e', e'' \}$.
Then,  it follows from the above proposition that, by means of the image of $\mr{Com} (e, e', e'')^{\ncong} \sqcup \mr{Com} (e, e', e'')^{\ncong}_{\mr{op}}$  via $\xi_{e, e', e''}^{\ncong}$,  one  can recognize     whether or not  the composite $b \circ a$ can be   defined  in $\overline{X}$ (i.e., the codomain of $a$ coincides with the domain of $b$).
\end{rema}
\SSP

Now let us 
 combine the previous two propositions (and Proposition \ref{L090}).
Let $G$, $G'$, $G''$, $e$, $e'$,  and $e''$ be as above.
Moreover, we will write
\begin{align}
\mr{Com}(e, e', e'') &:=\left\{ (u, u', u'') \in  \mr{Mor} (e, e')   \times \mr{Mor} (e', e'') \times  \mr{Mor}(e, e'') \, \Big| \, u'' = u' \circ u  \right\}, \\
\mr{Com} (e, e', e'')_\mr{op} &:=\left\{ (u, u', u'') \in  \mr{Mor}(e', e)   \times \mr{Mor} (e'', e') \times  \mr{Mor}(e'', e) \, \Big| \, u'' = u \circ u'  \right\}.  \notag
\end{align}
That is to say, $\mr{Com}(e, e', e'')$ and $\mr{Com} (e, e', e'')_\mr{op}$ are defined as the graphs of the composition law defining $\overline{X}$.

For each $e, e' \in \mr{Ob}(\overline{X})$, we write
\begin{align}
{^\dagger}\mr{Mor}\{ e, e' \} = \begin{cases} {^\dagger}\mr{Aut}(e) \sqcup {^\dagger}\mr{Mor}\{ e, e \}^\ncong & \text{if $e = e'$}, \\ {^\dagger}\mr{Mor}\{ e, e \}^\ncong  &  \text{if $e \neq e'$}.\end{cases}
\end{align}
If $e = e'$ (resp., $e \neq e'$) in $\overline{X}$, then we have $\mr{Mor}\{e, e' \} = \mr{Aut}(e) \sqcup \mr{Mor}\{ e, e \}^\ncong$ (resp., $\mr{Mor}\{e, e' \} = \mr{Mor}\{e, e' \}^\ncong$).
Hence,  the bijections $\xi_{(-)}^{\cong}$ (cf. (\ref{E1235})) and $\xi_{(-), (-)}^{\ncong}$ (cf. (\ref{E9898})) together 
yield a bijection
\begin{align} \label{E43928}
\xi_{e, e'} : \mr{Mor}\{ e, e' \} \isom  {^\dagger}\mr{Mor}\{e, e' \}.
\end{align}
Propositions \ref{L090}, \ref{L34}, and \ref{L07}  imply  the following  assertion.

\SSP
\bco \label{E67}
The image of $\mr{Com}(e, e', e'') \sqcup \mr{Com} (e, e', e'')_\mr{op}$  via the bijection
\begin{align}
\xi_{e, e', e''} := \xi_{e, e'} \times \xi_{e', e''} \times \xi_{e, e''} : & \ \mr{Mor} \{ e, e' \} \times  \mr{Mor} \{ e', e'' \} \times \mr{Mor} \{ e, e'' \} \\
& \  \isom {^\dagger}\mr{Mor} \{ e, e' \} \times  {^\dagger}\mr{Mor} \{ e', e'' \} \times {^\dagger}\mr{Mor} \{ e, e'' \} \notag
\end{align}
can be reconstructed category-theoretically from the data $(\mcC on_X, F^+_{G, e}, F^+_{G', e'}, F^+_{G'', e''})$, i.e., of a category and three distinct minimal objects in this category
  (cf. Propositions  \ref{L090} and  \ref{L06} for the category-theoretic reconstruction of ${^\dagger}\mr{Mor}\{ -, - \}$).
\eco

\LSP
\subsection{Proof of  Theorem \ref{E223}} \label{SS07}
We can  now prove Theorem \ref{E223} by applying the results obtained so far.
First, let us prove  the implication (a) $\Rightarrow$ (b).
To this end,  it suffices
 to show that $\mcS p_X \cong \mcS p_{X^\mr{op}}$, or equivalently $\mcC on_X \cong \mcC on_{X^\mr{op}}$ (cf. (\ref{E14})).
But this fact  can be  verified because the assignment from each constructive functor  $F :Y \migi X$  to $F^\mr{op} : Y^\mr{op} \migi X^\mr{op}$ yields an equivalence of categories $\delta_X : \mcC on_X \isom \mcC on_{X^\mr{op}}$.

Next, let us consider the inverse direction  (b) $\Rightarrow$ (a).
Suppose that there exists an equivalence of categories $\mcS p_X \isom \mcS p_{X'}$;
 this equivalence together with (\ref{E14}) determines  an equivalence of categories $\alpha : \mcC on_X \isom \mcC on_{X'}$.
Let us fix skeletons  of $X$ and  $X'$, which we denote by $\overline{X}$ and  $\overline{X}'$ respectively.
By  Proposition \ref{L01}, $\alpha$ induces, via $\xi_{\overline{X}}$ and $\xi_{\overline{X}'}$,   a bijection
\begin{align}
\overline{\alpha} :  \mr{Ob}(\overline{X}) \isom \mr{Ob}(\overline{X}').
\end{align}
It follows from Propositions \ref{L090} and \ref{L06} that, for each $e, e' \in \mr{Ob}(\overline{X})$, the equivalence  $\alpha$ induces, via $\xi_{e, e'}$ and $\xi_{\overline{\alpha}(e), \overline{\alpha} (e')}$,   a bijection
\begin{align}
\overline{\alpha}_{e, e'} : \mr{Mor} \{ e, e' \} \isom \mr{Mor} \{ \overline{\alpha} (e), \overline{\alpha} (e') \}.
\end{align}
Moreover, by applying Corollary \ref{E67} and using $\xi_{e, e', e''}$ and $\xi_{\overline{\alpha}(e), \overline{\alpha}(e'), \overline{\alpha}(e'')}$, we obtain a bijection
\begin{align}
\overline{\alpha}_{e, e', e''} : & \  \mr{Com}(e, e', e'') \sqcup \mr{Com} (e, e', e'')_\mr{op}   \\& \ \isom \mr{Com}(\overline{\alpha}(e), \overline{\alpha}(e'), \overline{\alpha}(e'')) \sqcup \mr{Com} (\overline{\alpha}(e), \overline{\alpha}(e'), \overline{\alpha} (e''))_\mr{op} \notag
\end{align}
for any $e, e', e'' \in \mr{Ob}(\overline{X})$.
 If  $\overline{X}$  
  has only one object $e$, then the bijections $\overline{\alpha}_e$, $\overline{\alpha}_{e, e}$, and $\overline{\alpha}_{e,e,e}$ show  that $\overline{X} \cong \overline{X}'$, which implies    $X \cong X'$.
 Hence,  it suffices to consider the case where 
 $\sharp \mr{Ob}(\overline{X}) \geq 2$.
 Since $X$ is assumed to be connected, there exists  
 a morphism $u_0 : e_0 \migi e'_0$ in $\overline{X}$ with $e_0 \neq e'_0$.

Here,   suppose that $\overline{\alpha}_{e_0, e'_0}(u) \in \mr{Mor}(\overline{\alpha} (e_0), \overline{\alpha}(e'_0))$.
 By the connectedness of  $X$ again,  one may verify (cf. Remark \ref{R340}) that $\overline{\alpha}_{e, e'}$ is restricted to a bijection 
  \begin{align}
  \overline{\alpha}_{e, e'}^\circledcirc : \mr{Mor}(e, e') \isom \mr{Mor}(\overline{\alpha}(e), \overline{\alpha}(e'))
  \end{align}
   for any $e, e' \in \mr{Ob}(\overline{X})$.
   Moreover, 
  $\overline{\alpha}_{e, e', e''}$ is restricted to a  bijection
  \begin{align}
  \overline{\alpha}_{e, e', e''}^\circledcirc &: \mr{Com}(e, e', e'')  \isom \mr{Com}(\overline{\alpha}(e), \overline{\alpha}(e'), \overline{\alpha}(e''))
  \end{align}
 for any $e, e', e'' \in \mr{Ob}(\overline{X})$.
 The   bijections $\overline{\alpha}$, $\overline{\alpha}_{e, e'}^\circledcirc$, and $\overline{\alpha}_{e, e', e''}^\circledcirc$ for various $e, e', e''$ together  yield an equivalence of categories
 $\overline{X} \isom \overline{X}'$.
 Thus, we conclude that  $X \cong  X'$.
 
 On the other hand, suppose that $\overline{\alpha}_{e_0, e'_0} (u_0) \in \mr{Mor} (\overline{\alpha} (e'_0), \overline{\alpha} (e_0))$.
 Write  $\alpha' := \delta_{X'} \circ \alpha$, where $\delta_{(-)}$ denotes the equivalence of categories constructed in the proof of (a) $\Rightarrow$ (b).
 Then, 
we have $\overline{\alpha}'_{e_0, e'_0} (u_0) \in \mr{Mor}_{\overline{X}'^{\mr{op}}}(\overline{\alpha} (e_0), \overline{\alpha} (e'_0))$.
 By
applying the above discussion to $\alpha'$,
we conclude that  $X \cong X'^{\mr{op}}$ (or equivalently, $X^\mr{op} \cong X'$).
This completes the proof of the implication (b) $\Rightarrow$ (a).
Thus,  we have  finished the proof of Theorem \ref{E223}.

\SSP
\begin{rema} \label{R45fg0}
The above proof 
 shows that, in order to reach the conclusion,  it suffices to assume the connectedness only for  either $X$ or $X'$.
In other words, the connectedness of $X$ can be characterized category-theoretically by  the  category $\mcS p_X$.
\end{rema} 

\SSP
\begin{rema} \label{R450}
The reconstruction carried out in Theorem \ref{E223} is reminiscent of ``Grothendieck conjecture"-type theorems in anabelian geometry. 
Anabelian geometry is, roughly speaking, an area of arithmetic geometry that  
 discusses  the issue of how much information concerning the geometry of certain arithmetic varieties  (e.g., hyperbolic curves  over a number field or a $p$-adic local field) is contained in the knowledge of the \'{e}tale fundamental groups, or equivalently, the categories of finite \'{e}tale coverings.
The classical point of view of anabelian geometry
 centers around a comparison between two arithmetic varieties  (or more generally, two   geometric objects of the same kind) via their \'{e}tale fundamental groups and it is  referred to as {\it bi-anabelian geometry}.
On the other hand, 
{\it mono-anabelian geometry}, being an alternative and  relatively new formulation,  centers around the task of establishing a 
group-theoretic algorithm whose input data consists of a single
abstract topological group isomorphic to the \'{e}tale fundamental group of a single
geometric object.
In particular, it requires us to reconstruct, unlike the bi-anabelian formulation, the desired data without any mention of some ``fixed reference model" copy of initial objects.
For basic references, we refer the reader to ~\cite{Hos}, ~\cite{Moc3}.

As explicitly verbalized in ~\cite{Son}, \S\,3,  constructive functors  on a category $X$ may be regarded, in some sense,  as categorical realizations of   coverings over $X$.
In fact,  if $X$ is a groupoid $G$, then the results of Proposition \ref{P014} enable us to consider the constructive functor $F_{G, e} : Y_{G, e} \migi G$ as if it were the universal covering of $G$.
Accordingly, we might say that  $\mcS p_X$ ($\cong \mcC on_X$ by (\ref{E14})) 
is like  the Galois   category consisting of finite  \'{e}tale coverings.
Notwithstanding the fact that Theorem \ref{E223} is stated in a completely bi-anabelian way,
the proof of the theorem actually furnishes a mono-anabelian algorithm that reconstructs $X$ from $\mcS p_X$. 
Here,  however, we will  omit the details of the formulation, as well as  the proof,  in that way.
\end{rema}
\SSP

\begin{rema} \label{R456}
C. L\'{e}vi-Strauss developed
 a structuralist theory of mythology which attempted to explain how seemingly fantastical and arbitrary tales could be so similar across cultures (cf. ~\cite{Lev1}, ~\cite{Lev2}). 
 Because he believed that there was no one authentic version of a myth, rather that they were all manifestations of the same language, he sought to find the fundamental units of myth, namely, the {\it mytheme}. 
The canonical formula of mythical transformation is an expression proposed in 1955 by  L\'{e}vi-Strauss in order to account for the abstract relations occurring between characters and their attributes in a myth understood as the collection of its variants. 
According to the canonical formula, a myth is reducible to an expression:
\begin{align} \label{EE39f}
f_x (a) : f_y (b) \cong f_x (b) : f_{a^{-1}}(y),
\end{align}
where each of these four arguments consist of a term variable ($a$ and $b$), and a function variable ($x$ and $y$).
This formula describes a structural relationship between a set of narrative terms and their transmutative relationships.
See ~\cite{Mar} for a reference on this formula.

Note that J. Morava tried to develop a truly mathematical argument for the canonical formula (cf. ~\cite{Mor1}, ~\cite{Mor2}).
He proposed to interpret it
 as the existence of an anti-isomorphism of the quaternion group.
On the other hand, as a mathematical study of myths in another direction, we might expect  that some kind of symmetry or structure on the whole category of structure species satisfying the condition expressed by (\ref{EE39f})  serves as a metaphorical explanation of some truth that L\'{e}vi-Strauss expected from myths. 
However, at the time of writing the present paper, the author does not have any effective ideas for the development of this argument.
\end{rema}

\vspace{10mm}
\section{Appendix: Category of categories over a category} 
\SSP

In this Appendix, we establish,   as  an analogy  of Theorem \ref{E223},
 the reconstruction of a category $X$ from  
 the category of categories over $X$, i.e., the category of  functors $F: Y \migi X$ (that are 
  not necessarily a structure species).
 The conclusion differs from the case of Theorem \ref{E223} in that we can reconstruct (up to a certain indeterminacy) not only the equivalence class of $X$ but also its isomorphism class.
 The following Proposition \ref{L11} (resp., \ref{L23}; resp., \ref{L24}) will be proved by using  a much simpler argument than  was used in the corresponding previous  assertion, i.e., Proposition  \ref{L01} (resp., \ref{L06}; resp., \ref{L07}).
So we will leave their  proofs to  the reader.

Let $X$ be a $\msU$-small category.
Denote by 
\begin{align}
\mcC at_X
\end{align} 
the category  defined as follows:
\begin{itemize}
\item 
The objects are functors
 of the form $F : Y \migi X$, where $Y$ is a $\msU$-small category.
\item
The morphisms from $F : Y \migi X$ to $F': Y' \migi X$ are functors $\Psi : Y \migi Y'$ satisfying $F = F' \circ \Psi$.
\end{itemize}
 The identity functor  $\mr{id}_X : X \isom X$ is a terminal object in this category.

For each object $e$ in $X$, we shall write $Y_e := \mcD is (\{ e \})$ and write 
\begin{align}
F_e : Y_e  \migiincl X
\end{align}
for the natural functor; this functor species an object in $\mcC at_X$.
We denote by
\begin{align}
{^\ddagger}\mr{Ob}(X) \left(\subseteq \mr{Ob}(\mcC at_X) \right)
\end{align}
the set of objects in $\mcC at_X$ of the form $F_e$ for some $e \in \mr{Ob}(X)$.
The assignment   $e \mapsto F_e$
 gives a bijection of sets
\begin{align} \label{Erex98}
\zeta_{X} : \mr{Ob}(X) \isom {^\ddagger}\mr{Ob}(X).
\end{align}

\SSP
\bpr \label{L11}
Let $F : Y \migi X$ be an object in $\mcC at_X$.
Then, $F$ is 
minimal in $\mcC at_X$
 if and only if $F$ is isomorphic to $F_e$ for some $e \in \mr{Ob}(X)$.
 In particular, the subset ${^\ddagger}\mr{Ob}(X)$ of $\mr{Ob} (\mcC at_X)$ can be reconstructed category-theoretically from the data $\mcC at_X$ (i.e., of a category).
\epr
\SSP

Next, let $e$ and  $e'$ be (possibly the same)  objects in $X$.
For each morphism  $v : e \migi e'$  in $X$, we shall  set
$Y_{v}$
to  be the subcategory of $X$ 
satisfying  $\mr{Ob}(Y_v) = \{ e, e' \}$ and $\mr{Mor}(Y_v) = \{ \mr{id}_e, \mr{id}_{e'}, v \}$.
Denote by  the natural inclusion by 
\begin{align}
 F_v :  Y_v \migiincl  X.
 \end{align}
The inclusions $Y_e \migiincl Y_v$ and $Y_{e'} \migiincl Y_v$ induce the morphism
 \begin{align}
 \Upsilon_v  : F_e \sqcup F_{e'} \migi F_v
 \end{align}
 in $\mcC at_X$.
 
We shall denote by
\begin{align}
\mcD_{e, e'}
\end{align}
the category defined as follows:
\begin{itemize}
\item
The objects in $\mcD_{e, e'}$ are pair $(F, \Upsilon)$ consisting of an object $F : Y \migi X$ in $\mcC at_X$ and a non-invertible monomorphism $\Upsilon : F_e \sqcup  F_{e'} \migiincl F$.
\item
The morphisms from $(F, \Upsilon)$ to $(F', \Upsilon')$ are morphisms $\Psi : F \migi F'$ satisfying $\Upsilon' = \Psi \circ \Upsilon$.
\end{itemize}
For each element $v \in \mr{Mor}(e, e')$,
the pair $(F_v, \Upsilon_v)$ introduced above specifies an object in $\mcD_{e, e'}$.
The assignment $(F, \Upsilon) \mapsto \Upsilon$ defines a functor
\begin{align} \label{E449}
\mcD_{e, e'} \migi \mcC at_X.
\end{align}
It follows from Proposition \ref{L11} that the category $\mcD_{e, e'}$ together  with the functor (\ref{E449}) can be reconstructed category-theoretically from the data $(\mcC at_X, F_e, F_{e'})$ (i.e.,  of a category and two minimal objects in this category).
Also, as in the case of ``$\mcC_{e, e'}$" introduced in \S\,\ref{SS58},  there exists a natural identification $\mcD_{e, e'} = \mcD_{e', e}$, by which we will not distinguish between $\mcD_{e, e'}$ and $\mcD_{e', e}$.
In particular, this  identification allows us to consider $(F_w, \Upsilon_w)$'s for various $w \in \mr{Mor}(e', e)$ as objects in $\mcD_{e, e'}$. 
We denote by
\begin{align}
{^\ddagger}\mr{Mor} \{e, e' \}
\end{align}
the set of morphisms in $\mcC at_X$ of the form $\Upsilon_v$ for  some  $v \in \mr{Mor}\{e, e' \}$.
The assignment $v \mapsto \Upsilon_v$ defines a bijection of sets
\begin{align}
\zeta_{e, e'} : \mr{Mor}\{ e, e' \} \isom  {^\ddagger}\mr{Mor}\{ e, e' \}.
\end{align}

\SSP
\bpr \label{L23}
Let $(F, \Upsilon)$ be an object in $\mcC at_X$.
The following two conditions are equivalent to each other:
\begin{itemize}
\item
$(F, \Upsilon)$ is minimal in $\mcD_{e, e'}$ and there is no triple of morphisms
\begin{align}
(\Upsilon_0: F_e \migi F_0, \Upsilon'_0: F_{e'} \migi F'_0, \Psi : F_0 \sqcup F'_0 \migi F)
\end{align}
in $\mcC at_X$ such that $\Psi$ is an isomorphism and satisfies the equality $\Upsilon = \Psi \circ (\Upsilon_0 \sqcup \Upsilon'_0)$.
\item
$(F, \Upsilon)$ is isomorphic to $(F_v, \Upsilon_v)$ for some $v \in \mr{Mor}\{ e, e' \}$.
\end{itemize}
In particular, the subset ${^\ddagger}\mr{Mor}\{ e, e' \}$ of $\mr{Mor}(\mcC at_{X})$ can be reconstructed category-theoretically from the data $(\mcC at_X, F_e, F_{e'})$ (i.e., of a category and two minimal objects in this category).
\epr
\SSP

Given  three (possibly the same) objects $e, e', e''$ in $X$, 
we shall write
\begin{align} \label{E0090}
\mr{Com}(e, e', e'')^{\mcC at} &:=\left\{ (u, u', u'') \in  \mr{Mor}(e, e')\times \mr{Mor} (e', e'') \times  \mr{Mor} (e, e'') \, \Big| \, u'' = u' \circ u  \right\},\\
\mr{Com}(e, e', e'')_\mr{op}^{\mcC at} &:=\left\{ (u, u', u'') \in   \mr{Mor}(e', e) \times  \mr{Mor} (e'', e') \times   \mr{Mor} (e'', e) \, \Big| \, u'' = u \circ u'  \right\}.
 \notag
\end{align}
Then, we have the following assertion.

\SSP
\bpr \label{L24}
Let us choose $u \in \mr{Mor}\{e, e' \}$, $u' \in \mr{Mor} \{ e', e'' \}$, and  $u'' \in \mr{Mor}\{ e, e'' \}$.
Then, the  following two conditions are equivalent to each other:
\begin{itemize}
\item
Either $u \circ u'$ or $u' \circ u$ can be defined and one of the equalities $u'' = u \circ u'$, $u'' = u' \circ u$ holds.
\item
The colimit  $F_{\mr{limit}}^{\mcC at}$ of the diagram
\begin{align} \label{E446}
\vcenter{\xymatrix@C=46pt@R=36pt{
F_e \ar[r]^-{\mr{inclusion}} \ar[d]_-{\mr{inclusion}} & F_{u'} \\
F_u &
}}
\end{align}
exists in $\mcC at_X$, and there exists a monomorphism $F_{u''} \migiincl F_{\mr{limit}}^{\mcC at}$ in $\mcC at_X$  via which the morphism $\Upsilon_{u''} : F_e \sqcup F_{e''} \migiincl F_{u''}$ is compatible with the morphism $F_e \sqcup F_{e''} \migiincl F_{\mr{limit}}^{\mcC at}$ arising naturally from $\Upsilon_u$ and $\Upsilon_{u'}$.
\end{itemize}
In particular, the image of $\mr{Com}(e, e', e'')^{\mcC at} \sqcup \mr{Com}(e, e', e'')^{\mcC at}_{\mr{op}}$ via the bijection
\begin{align}
\zeta_{e, e', e''} := \zeta_{e, e'} \times \zeta_{e', e''} \times \zeta_{e, e''} : & \ \mr{Mor}\{e, e' \} \times  \mr{Mor}\{e', e'' \} \times \mr{Mor}\{e, e'' \} \\
& \ \isom   {^\ddagger}\mr{Mor}\{e, e' \} \times  {^\ddagger}\mr{Mor}\{e', e'' \} \times {^\ddagger}\mr{Mor}\{e, e'' \} \notag
\end{align}
can be reconstructed category-theoretically from the data $(\mcC at_X, F_e, F_{e'}, F_{e''})$ (i.e., of a category and three minimal objects in this category).
\epr
\SSP

The previous three propositions enable us to prove  the following assertion.
(The proof is entirely similar to the proof of Theorem \ref{E223}.)

\SSP
\bt \label{T05}
Let $X$ and $X$ be connected $\msU$-small categories.
Then, the following conditions are equivalent to each other:
\begin{itemize}
\item[(a)]
$X \stackrel{\mr{isom}}{\cong} X'$ or $X^\mr{op}  \stackrel{\mr{isom}}{\cong} X'$.
\item[(b)]
$\mcC at_X \cong \mcC at_{X'}$.
\end{itemize}
\et
\begin{proof}
The implication (a) $\Rightarrow$ (b) follows immediately from 
 the existence of an equivalence of categories $\mcC at_X \isom \mcC at_{X^\mr{op}}$ obtained  by assigning
 $F \mapsto F^\mr{op}$ for each $F \in \mcC at_X$.

Next, we shall consider the  inverse direction (b) $\Rightarrow$ (a).
Suppose that there exists an equivalence of categories $\beta : \mcC at_X \isom \mcC at_{X'}$.
By Proposition \ref{L11}, $\beta$ induces, via $\zeta_{X}$ and $\zeta_{X'}$,  a bijection
\begin{align}
\overline{\beta} : \mr{Ob}(X) \isom \mr{Ob}(X').
\end{align}
It follows from Proposition \ref{L23} that, for each $e, e' \in \mr{Ob}(X)$, the equivalence $\beta$ induces, via $\zeta_{e, e'}$ and $\zeta_{\overline{\beta}(e), \overline{\beta}(e')}$, 
 a bijection  of sets 
\begin{align}
\overline{\beta}_{e, e'} : \mr{Mor}\{e, e'\}  \isom \mr{Mor} \{\overline{\beta} (e), \overline{\beta} (e')\}. 
\end{align}
Moreover, by applying Proposition \ref{L24} and using the bijections $\zeta_{e, e', e''}$ and $\zeta_{\overline{\beta} (e), \overline{\beta} (e'), \overline{\beta} (e'')}$,
we obtain a bijection 
\begin{align}
\overline{\beta}_{e, e', e''} : & \ \mr{Com}(e, e', e'')^{\mcC at} \sqcup   \mr{Com}(e, e', e'')^{\mcC at}_{\mr{op}}  \\
& \ \isom  \mr{Com}(\overline{\beta} (e), \overline{\beta} (e'), \overline{\beta} (e''))^{\mcC at} \sqcup   \mr{Com}(\overline{\beta} (e), \overline{\beta} (e'), \overline{\beta} (e''))^{\mcC at}_{\mr{op}} \notag
\end{align}
  for any $e, e', e'' \in \mr{Ob}(X)$.
  If $X$ has only one object $e$, then the bijections $\overline{\beta}$,  $\overline{\beta}_{e, e}$, and $\overline{\beta}_{e, e,e}$ show   that 
   $X \stackrel{\mr{isom}}{\cong} X'$.
  Hence, it suffices to consider the case where $\sharp \mr{Ob}(X) \geq 2$.
  Since $X$ is assumed to be connected, 
there exists  a morphism $u_0 : e_0 \migi e'_0$ in $X$ with $e_0 \neq e'_0$.
  
  Here, suppose that $\overline{\beta}_{e_0, e'_0} (u_0) \in \mr{Mor}(\overline{\beta}(e_0), \overline{\beta}(e'_0))$.
  From the connectedness of  $X$ again, 
  one may verify (for the same reason as stated in  the comment on  Remark \ref{R340}) that 
  $\overline{\beta}_{e, e'}$ is restricted to a bijection
  \begin{align}
  \overline{\beta}^\circledcirc_{e, e'} : \mr{Mor} (e, e') \isom \mr{Mor}(\overline{\beta}(e), \overline{\beta}(e'))
  \end{align}
  for any $e, e' \in \mr{Ob}(X)$.
 Moreover, $\overline{\beta}_{e, e', e''}$ is restricted to a bijection
 \begin{align}
 \overline{\beta}_{e, e', e''} : \mr{Com}(e, e', e'')^{\mcC at} \isom 
 \mr{Com}(\overline{\beta} (e), \overline{\beta} (e'), \overline{\beta} (e''))^{\mcC at}
 \end{align} 
 for any $e, e', e'' \in \mr{Ob}(X)$.
 The bijections $\overline{\beta}$, $\overline{\beta}^\circledcirc_{e, e'}$, and $\overline{\beta}^\circledcirc_{e, e', e''}$ for various $e, e', e'' \in \mr{Ob}(X)$ yield an isomorphism of categories $X \isom X'$.
 
 On the other hand,  if  $\overline{\beta}_{e_0, e'_0} (u_0) \in \mr{Mor}(\overline{\beta}(e'_0), \overline{\beta}(e_0))$,
then the problem, as in the proof of Theorem \ref{E223},  reduces   to the previous case by using the equivalence of categories $\mcC at_{X'} \isom \mcC at_{X'^{\mr{op}}}$ given by $F \mapsto F^{\mr{op}}$.
 At any rate, we conclude that $X \stackrel{\mr{isom}}{\cong} X'$ or $X^{\mr{op}} \stackrel{\mr{isom}}{\cong} X'$, as desired. 
 This completes the implication (b) $\Rightarrow$ (a).
  \end{proof}

\end{document}